\documentclass[11pt]{article}
\usepackage{amsmath,amssymb,bbm}

\makeatletter
\newfont{\footsc}{cmcsc10 at 8truept}
\newfont{\footbf}{cmbx10 at 8truept}
\newfont{\footrm}{cmr10 at 10truept}
\newtheorem{Definition}{Definition}

\newtheorem{Lemma}{Lemma}

\newcommand{\tmop}[1]{\ensuremath{\operatorname{#1}}}

\newtheorem{Theorem}{Theorem}

\newenvironment{proof}{
  \noindent\textbf{Proof}\ }{\hspace*{\fill}
  \begin{math}\Box\end{math}\medskip}
\begin{document}
\title{Random walks on hypergroup of conics in finite fields}

\author{Le Anh Vinh\\
\small School of Mathematics,\\[-0.8ex]
\small University of New South Wales\\[-0.8ex]
\small Sydney 2052 NSW\\[-0.8ex]
\small \texttt{vinh@maths.unsw.edu.au}}
\date{\empty}
\maketitle

\begin{abstract}
In this paper we study random walks on the hypergroup of conics in finite fields. We investigate the behavior of random walks on this hypergroup, the equilibrium distribution and the mixing times. We use the coupling method to show that the mixing time of random walks on hypergroup of conics is only linear.
\end{abstract}
\begin{center}

\small Keywoords: random walks, hypergroups, finite fields.
\end{center}

\section{Introduction}

Throughout this paper, $q$ is an odd prime power, $F_q = GF(q)$ is the finite field with $q$ elements, $F_q^\ast$ is the multiplicative group of the non-zero elements of $F_q$ and $a, b,c$ are any three fixed numbers in $F_q^\ast$ such that $ab = c^2$.

\begin{Definition}
  The \textbf{weighted-quadrance} $Q^{a,b} ( A_1, A_2 )$ between the points 
  $A_1 = ( x_1, y_1)$, and $A_2 = ( x_2, y_2 )$ is the number
  \[ Q^{a,b} ( A_1, A_2 ) := a( x_2 - x_1 )^2 + b( y_2 - y_1 )^2 . \]
\end{Definition} 

For $a = b = 1$, we have the standard definition of quadrance which is introduced by Wildberger. The important of this notation is developed in \cite{7}.

\begin{Definition}
  A \textbf{weighted-circle} $C_k^{a,b}(A_0)$ in a finite field $F_q$ with center $A_0 \in F_q \times F_q$ and quadrance $k \in F_q$ is set of all points $X$in $F_q \times F_q$ such that
  \[ Q^{a,b}(A_0, X) = k. \]
\end{Definition}

Note that this notation generalizes those of ellipse and hyperbola in the classical setting and those of circle, quadrola and grammola in Wildberger's setting. 

We define $C_i^{a,b}$ the weighted-circle centered at the origin and quadrance $i \in F_q$. Let $N_i^{a,b}$ be the number of solutions of the equation $ax^2+by^2 = i$ in the field $F_q$. Then $N_i^{a,b}$ is the number of points in $C_i^{a,b}$. Hence, we have a partition of the finite field space $F_q^2$ into $q$ set of points $\{ C_i^{a,b} \}_{i\in F_q}$. If we start from $O=(0,0)$, take a random step by translating by
an element of $C_i^{a,b}$, and then take another random step by translating by an
element of $C_j^{a,b}$, the final point will be an element of $C_k^{a,b}$ for some
$k$. Counting over all possible such combinations, there are $N_{ij}^{k} $
ways to reach to a point of $C_{k}^{a,b}$ by using steps from $C_{i}^{a,b}$ then $C_{j}^{a,b}$ randomly. 
We can write this relation as
\begin{equation*}C_i^{a,b}C_j^{a,b}=\sum_{k \in F_q}N_{ij}^{k}C_k^{a,b},\end{equation*}
where $N_{ij}^{k}$ are non-negative integers. 

Let $n_{ij}^{k}=\frac{N_{ij}^{k}
}{|g_{i}\Vert g_{j}|}$ then this relation can be written as distribution form
\begin{equation}\label{eq1} C_i^{a,b} C_j^{a,b} = \sum_{k \in F_q} n_{i j}^{k} C_k^{a,b},\end{equation}
where $n_{i j}^{k} \geqslant0$ and $\sum_{k} n_{i j}^{k} = 1$ for any $i, j$.

We recall the formal definition of (general) hypergroup (see \cite{6}).

\begin{Definition}
  A (finite) general hypergroup is a pair $(\mathcal{K},\mathcal{A})$ where
  $\mathcal{A}$ is a *-algebra with unit $c_0$ over $\mathbbm{C}$ and
  $\mathcal{K}= \{ c_0, c_1, \ldots, c_n \}$ is a subset of $\mathcal{A}$
  satisfying
  \begin{enumerate}
    \item $\mathcal{K}$ is a basis of $\mathcal{A}$    
    \item $\mathcal{K}^{\ast} =\mathcal{K}$    
    \item The struture constants $n_{i j}^k \in \mathbbm{C}$ defined by    
     \[ c_i c_j = \sum_k n_{i j}^k c_k \]
    satisfy the conditions    
     \begin{align*}
     c_i^{\ast} &= c_j \Leftrightarrow n_{i j}^0 > 0,\\
     c_i^{\ast} &= c_j \Leftrightarrow n_{i j}^0 > 0.
     \end{align*}      
  \end{enumerate}
\end{Definition}

$\mathcal{K}$ is called \textit{hermitian} if $c_i^{\ast} = c_i$ for all $i$,
\textit{commutative} if $c_i c_j = c_j c_i$ for all $i, j$, \textit{real} if $n_{i j}^k$ $\in R$ for all $i, j, k$, \textit{positive} if $n_{i j}^k \geqslant 0$ for all $i, j, k$ and \textit{normalized} if $\sum_k n_{i j}^k = 1$ for all $i, j$. A generalized hypergroup which is both positive and normalized will be called a hypergroup. There are board examples and applications of (generalized) hypergroups which can be found in \cite{6}.

In \cite{3}, we studied the random walk over hypergroup of circles (i.e. $a = b = 1$) in a finite
field of prime order $p = 4l + 3$ using comparison of Dirichlet Forms and geometric bound of
eigenvalues for Markov chains. In this paper, we will study the random walk over hypergroup of weighted-circles (more general conics) in arbitrary finite field using the coupling method. In Section 2, we will show that the set $\mathcal{C} = \{ C_i^{a,b}\}_{i \in F_q}$ with the relation \ref{eq1} is a hypergroup. The structure of this hypergroup will also be given. In Section 3, we will study the random walk over the hypergroup of weighted-circles. The main result of this paper is that the random walk over hypergroup of weighted-circle tends to the stationary distribution in a linear time with respect to the size of the hypergroup.

\section{Hypergroup of weighted-circles}

\subsection{Some Lemmas}

Recall that a (multiplicative) character of $F_q$ is a homomorphism from
$F_q^{\ast}$, the multiplicative group of the non-zero elements of $F_q$, to the
multiplicative group of complex numbers with modulus $1$. The identically $1$
function is the principal character of $F_q$ and is denoted $\chi_0$. Since
$x^{q - 1} = 1$ for every $x \in F_q^{\ast}$ we have $\chi^{q - 1} = \chi_0$
for every character $\chi$. A character $\chi$ is of order $d$ if $\chi^d =
\chi_0$ and $d$ is the smallest positive integer with this property. By
convention, we extend a character $\chi$ to the whole of $F_q$ by putting
$\chi ( 0 ) = 0$. The quadratic (residue) character is defined by $\chi ( x )
= x^{( q - 1 ) / 2}$. Equivalently, $\chi$ is $1$ on square, $0$ at $0$ and $-
1$ otherwise.

The following lemma give us the number of points in any weighted-circle in $F_q^2$.
\begin{Lemma} \label{number points}
  If $i \in F_q^\ast$ then \[N_i^{a,b} = q - (-1)^{(q-1)/2}.\]
\end{Lemma}
\begin{proof} 
From Proposition 8.6.1 and Theorem 5 in \cite[p. 101-104]{number theory}, we know that
\[
	N_i^{a,b} = q + (-1)^{(q-1)/2} \chi(a^{-1}) \chi(b^{-1}).
\]
But $a b = c^2$ so $\chi(a^{-1}) \chi(b^{-1}) = 1$. This concludes the proof of the lemma.
\end{proof}

The following lemma gives us the number of intersections between any two weighted-circles in $F_q^2$.

\begin{Lemma} \label{khac khong}
Let $i, j \ne 0$ in $F_q$ and let $X, Y$ be two distinct points in $F_q^2$ such that $Q^{a,b}(X,Y) = k \ne 0 $. Then $\left|C_i^{a,b}(X) \cap C_j^{a,b}(Y)\right|$ only depends on $i, j$ and $k$. Precisely,  let $f ( i, j, k ) = i j - (i - j - k )^2/4$. Then the number of intersection points is $p_{i j}^k$, where
  \begin{equation}\label{eq2}
  p_{ij}^k =
   \begin{cases}
    0  & \text{if}\; \; f(i,j,k)\; \text{is non-square},\\
    1  & \text{if}\; \; f(i,j,k) = 0,\\
    2 &\text{if}\; \; f(i,j,k)\; \text{is square}.
   \end{cases}
  \end{equation}

\end{Lemma}

\begin{proof}
Suppose that $X = ( m, n )$ and $Y = ( m + x, n + y )$ for some $m, n, x, y \in F_q$ then $ax^2 + by^2 = k$. Suppose that $Z \in C_i^{a,b}(X) \cap C_j^{a,b}(Y)$ where $Z = (m + x + u, n + y + v)$ for some $u, v \in F_q$. Then we have $au^2 + bv^2 = j$ and $a( x + u )^2 + b( y + v^2 ) = i$. This implies that $ax u + by v = (i - j - k)/2$.
But we have $( ax u + by v )^2 + ( cx v - cy u )^2 = ( ax^2 + by^2 ) ( au^2 + bv^2)$ so
  \begin{align*} ( cx v - cy u )^2 &= k j - ( i - j - k )^2/4 \\
                                   &= i j - ( k - i - j )^2/4 = f ( i, j, k ) . \end{align*}
  If $f ( i, j, k )$ is non-square number in $F_q$ then it is clear that there does not exist
  such $x, y, u, v$, or $p_{i j}^k = 0$. Otherwise, let $\alpha = ( i - j
  - k ) / 2$ and $f ( i, j, k ) = \beta^2$ for $0 \leqslant \beta \leqslant (
  p + 1 ) / 2$ then
  \[
    cx v - cy u = \pm \beta,\qquad ax u + by v = \alpha.
  \]
  Solving for $( u, v )$ with respect to $( x, y )$ we have
  \[
    u = \frac{\alpha cx \mp \beta by}{ck},\qquad  v = \frac{\alpha cy \pm \beta ax}{ck}. 
  \]
  If $\beta = 0$ then we have only one $( u, v )$ for each $( x, y )$, but if
  $\beta \neq 0$ then we have two pairs $( u, v )$. This implies (\ref{eq2}), completing the proof.
\end{proof}

\subsection{The first case}

Suppose that $q \equiv 3$ (mod $4$). From Lemma \ref{number points}, we have
  \begin{equation*}
  N_{i}^{a,b} =
   \begin{cases}
    q+1  & \text{if}\; \; i\in F_q^\ast,\\
    1 &\text{otherwise}.
   \end{cases}
  \end{equation*}
  
Since $N_0^{a,b}=1$ so $Q^{a,b}(X,Y)=0$ if and only if $X \equiv Y$. Hence for any $j \in F_q$ then $C_0^{a,b}C_j^{a,b}=C_j^{a,b}$. The following theorem details the coefficients in (\ref{eq1}).
\begin{Theorem} \label{hyper 3} Suppose that $q \equiv 3$ (mod $4$), $i,j \in F_p^\ast$ and $k \in F_q$. Then $n_{ij}^{k}$ (in (\ref{eq1})) only depend on $f(i,j,k)$.
Precisely, we have
  \begin{equation*}
  n_{ij}^k =
   \begin{cases}
    0  & \text{if}\; \; f(i,j,k)\; \text{is non-square},\\
    1/(q+1)  & \text{if}\; \; f(i,j,k) = 0,\\
    2/(q+1) &\text{if}\; \; f(i,j,k)\; \text{is square}.
   \end{cases}
  \end{equation*}
\end{Theorem}

\begin{proof}
There are there cases. 

\begin{enumerate}
    \item Suppose that $f(i,j,k)$ is non-square. From Lemma \ref{khac khong}, for any $( x, y )$ in $C_i^{a,b}$, there 
    does not exists $( u, v )$ in $C_j^{a,b}$ such that if we go by $( x, y )$ followed by $( u, v )$, the destination 
    is a point in $C_k^{a,b}$. Hence $N_{i,j}^k = 0$ and $n_{i,j}^k = 0$.
    \item Suppose that $f(i,j,k) = 0$. From Lemma \ref{khac khong}, for any $( x, y )$ in $C_i^{a,b}$, there 
    exists a unique $( u, v )$ in $C_j^{a,b}$ such that if we go by $( x, y )$ followed by $( u, v )$, the destination 
    is a point in $C_k^{a,b}$. Hence $N_{i,j}^k = |C_i^{a,b} | = p + 1$ and $n_{i,j}^k = 1/(p+1)$.
    \item Suppose that $f(i,j,k)$ is square. From Lemma \ref{khac khong}, for any $( x, y )$ in $C_i^{a,b}$, there 
    exists two points $( u, v )$ in $C_j^{a,b}$ such that if we go by $( x, y )$ followed by $( u, v )$, the
    destination is a point in $C_k^{a,b}$. Hence $N_{i,j}^k = 2|C_i^{a,b} | = 2(p + 1)$ and $n_{i,j}^k = 2/(p+1)$.
  \end{enumerate}
  This concludes the proof of the theorem.
\end{proof}

Now, we show that the set $C = \{ C_i^{a,b}\}_{i \in F_q}$ with the
relation (\ref{eq1}) is a hypergroup. It is clear that $n_{i j}^k \geqslant
0$, and $\sum_{k\in F_q} n_{i j}^k = 1$ for any $i, j \in F_q$. From Theorem \ref{hyper 3}, $n_{i j}^0 \neq 0$ if and only if $f(i,j,0) =  - ( i - j )^2/4$ is square. But $q \equiv 3$ (mod 4) so $-1$ is not a square in $F_q$. Hence $n_{i j}^0 \neq 0$ if and only if $i = j$. Let $(C_i^{a,b})^{\ast} = C_i^{a,b}$ then $C$ is a hermitian commutative hypergroup (note that, $n_{i j}^k$ is symmetric with respect to $i, j$ and $k$ so $C$ is commutative).

\subsection{The second case}
Suppose that $q \equiv 1$ (mod $4$). From Lemma \ref{number points}, we have
  \begin{equation*}
  N_{i}^{a,b} =
   \begin{cases}
    q-1  & \text{if}\; \; i\in F_q^\ast,\\
    2q-1 &\text{otherwise}.
   \end{cases}
  \end{equation*}

This case, however, is harder since the null-circle $C_0^{a,b}$ contains more than one point and the set $C = \{C_i^{a,b}\}_{i\in F_q}$ turns out to be not a (hermitian) hypergroup. To resolve this difficulty, we need to redefine the null-circle $C_0^{a,b}$. We divide the null-circle into two parts 
\begin{align*}
 C_0^{a,b}(X) &= \{X\},\\
 C_q^{a,b}(X) &= \{ Y \in F_q^2 \;|\; Q^{a,b}(Y,X)=0, Y \ne X\}.
\end{align*}

We define $F_q^{+} = F_q \cup \{q\}$. The following theorem is similar to Theorem \ref{hyper 3}. The proof of this theorem is omitted as it is lengthly and repeated.
\begin{Theorem} \label{hyper 1}
  Suppose that $q \equiv 1$ (mod $4$) and $i,j, k \in F_q$. 
  \begin{enumerate}
    \item Suppose that $i = 0$. Then
    \begin{equation*}
    n_{0,j}^{k} =
   \begin{cases}
    1  & \text{if}\; \; j = k,\\
    0 &\text{otherwise}.
   \end{cases}
  \end{equation*}
    
    \item Suppose that $i,j \in F_q^\ast$.
    \begin{enumerate}
      \item If $k \in F_q^\ast$ then
    \begin{equation*}
  	n_{ij}^k =
   	\begin{cases}
    	0  & \text{if}\; \; f(i,j,k)\; \text{is non-square},\\
    	1/(q-1)  & \text{if}\; \; f(i,j,k) = 0,\\
    	2/(q-1) &\text{if}\; \; f(i,j,k)\; \text{is square}.
   	\end{cases}
  	\end{equation*}
     \item If $k = 0$ then
   \begin{equation*}
    n_{i,j}^{0} =
   \begin{cases}
    1/(q-1)  & \text{if}\; \; i = j,\\
    0 &\text{otherwise}.
   \end{cases}
  \end{equation*}
  \item If $k = q$ then
  \begin{equation*}
    n_{i,j}^{q} =
   \begin{cases}
    0  & \text{if}\; \; i = j,\\
    2/(q-1) &\text{otherwise}.
   \end{cases}
  \end{equation*}
  
    \end{enumerate}
\item Suppose that $i = q$.
    \begin{enumerate}
      \item If $j \in F_q^\ast$ then
        \begin{equation*}
    n_{q,j}^{k} =
   \begin{cases}
    1/(q-1)  & \text{if}\; \; j \ne k,\\
    0 &\text{otherwise}.
   \end{cases}
  \end{equation*}
      \item If $j = q$ then
        \begin{equation*}
    n_{q,q}^{k} =
   \begin{cases}
    1/2(q-1)  & \text{if}\; \; k \in F_q,\\
    (q-2)/2(q-1) &\text{otherwise}.
   \end{cases}
  \end{equation*}      
    \end{enumerate}
  \end{enumerate}
\end{Theorem}

From Theorem \ref{hyper 1}, it is clearly that $n_{i j}^0 > 0$ if and only if $i = j$. Hence the set $C =\{ C_i \}_{i\in F_q^{+}}$ with the random walk multiplication is a hermitian commutative hypergroup.

\section{Random walks over hypergroup of conics}
\subsection{Preliminary}
In this section, we will consider the random walk by $C_1^{a,b}$; that is we choose all steps from the weighted-circle $C_1^{a,b}$. For convenient, we drop the superscripts $a,b$ of weighted-circle and call this random walk $C_1$. This random walk has the Markov kernel $C_1(C_i,C_j)=n_{i,1}^{j}$ for all $i, j\in F_q$ (or $F_q^{+}$). In general, at $n^{\text{th}}$ step we have the relation
\[ C_1^n = \sum_{j\in F_q} \alpha_{n, j} C_j \]
where $\alpha_{n, j} \geqslant 0$ for $j\in F_q$ and $\sum_{j \in F_q} \alpha_{n,j}$ = 1.

Let $K$ be a Markov kernel. The probability $\pi$ is invariant or stationary
for $K$ if $\pi K = \pi$. A Markov kernel $K$ is irreducible if for any two
states $x, y$ there exists an integer $n = n ( x, y )$ such that $K^n ( x, y )
> 0$. A state $x$ is called aperiodic if $K^n ( x, x ) > 0$ for all
sufficiently large $n$. If $K$ is irreducible and has an aperiodic state then
all states are aperiodic and $K$ is \textit{erogodic}.

The following definition gives us the total variation distance between two
probability measures.

\begin{Definition}
  Let $\mu, \nu$ be two probability measures on the set $X$. The total
  variation distance is defined by  
  \begin{align*} d_{\tmop{TV}} ( \mu, \nu ) &= \max_{_{A
  \subset X}} | \mu ( A ) - \nu ( A ) | \\
  &= \frac{1}{2} \sum_{x \in X} | \mu ( x
  ) - \nu ( x ) |. \end{align*}
\end{Definition}

Ergodic Markov chains are useful algorithmic tools in which, regardless of
their initial state, they eventually reach a unique stationary distribution.
The following theorem, originally proved by Doeblin, details the essential
property of ergodic Markov chains.

\begin{Theorem}\label{ergodic}
  Let $K$ be any ergodic Markov kernel on a finite state space $X$ then $K$
  admits a unique stationary distribution $\pi$ such that
 \[
  \forall x, y \in X, \lim_{t \rightarrow
  \infty} K_t ( x, y ) = \pi ( y ). \]

\end{Theorem}

Let $K$ be a Markov chain on the set $X$ with the stationary distribution $\pi$. We define the \textit{mixing time} $\tau_p ( \varepsilon )$ as the time until the chain is within variation distance $\varepsilon$ to the distribution $\pi$, start from the worst initial state. We give a formal definition for this concept.

\begin{Definition} The mixing time $\tau_p(\varepsilon)$ is defined by
  \[\tau_p ( \varepsilon ) = \max_{x\in X} \min \{ t :
  d_{\tmop{TV}} ( K^t ( x , . ), \pi ) \leqslant \varepsilon \}.\]
\end{Definition}

We can fix $\varepsilon$ as any small constant. A popular choice is to set
$\varepsilon = 1/2 e$. We then boost to arbitrary small variation
distance by the following lemma.

\begin{Lemma}(\cite{1})
  $\tau_p ( \varepsilon ) \leqslant \tau_p ( 1/2 e ) \ln ( 1/e )$.
\end{Lemma}

\begin{Lemma}\label{lap}(\cite{1})
  Let $K$ be a Markov chain and $\pi$ be a probability distribtution on $X$. Suppose that there exists an integer $m$ 
  and a constant $c > 0$ such that for  all $x, y \in X$, $K^m ( x, y ) \geqslant c \pi ( y )$. Then
  $d_{\tmop{TV}} ( K^{m n} ( x, . ), \pi ) \leqslant ( 1 - c )^n$ for all
  integer $n$ and $x \in X$.
\end{Lemma}

\subsection{Main results}
\subsubsection{The first case}
Suppose that $q \equiv 3$ (mod $4$). We have the following lemma.

\begin{Lemma}\label{support 3}
	If $i, j \ne 0$ then there exists $k$ such that $n_{i1}^k, n_{k1}^j >0$.
\end{Lemma}

\begin{proof}
  We start from a point in $C_i$, go a step from $C_1$ then we have $|C_i \| C_1|
  = ( q + 1 )^2$ possible steps. From Theorem \ref{hyper 3}, there is no more than $2 ( q
  + 1 )$ steps that can reach the weighted-circle. Thus, we reach
  at least $( q + 1 )^2/2 ( q + 1 ) = (q + 1)/2$ weighted-circles.  Applying the same argument, start from a point in
  $C_l$, go a step from $C_1$ then we reach at least $(q + 1)/2$ weighted-circles.
  Since we have only $q$ weighted-circles, by the pigeonhole principle, there exists a weighted-circle $C_k$
  which is reachable from both directions. The Lemma follows.
\end{proof}

From Lemma \ref{support 3}, we can show that the random walk $C_1$ is erogodic.  

\begin{Lemma}\label{important 3}         
  Let $q \equiv 3$ (mod $4$) and $\pi$ be the distribution on the set $C = \{ C_i\}_{i \in F_q}$ with $\pi(C_0)=1/q^2$ and $\pi(C_i) = (q+1)/q^2$ for all $i \in F_q^\ast$. Then
  \[
  C_1^4 ( C_i, C_j ) \geqslant \frac{q^2 ( q - 1 )}{( p +1
  )^4} \pi ( C_j ) \]
  for all $i, j \in F_q$.
\end{Lemma}
\begin{proof}
  There are four cases.
  \begin{enumerate}
  \item Suppose that $i = j = 0$. There are $|C_1 |^4 = ( q + 1 )^4$ possible ways to go by 4 steps. We
  first go by any 2 steps. In the last two, we just go backward then it is clear that we go
  back to the starting point. Therefore, at least $|C_1 |^2 = ( q + 1 )^2$
  ways to go from $C_0$ to $C_0 $. It implies that
 \[
  C_1^4 ( C_0, C_0 ) \geqslant \frac{(
  q + 1 )^2}{( q + 1 )^4} > \frac{q^{2} ( q - 1 )}{( q+1 )^4} \pi ( C_0 ).
 \]  
  
  \item Suppose that $i = 0, j \neq 0$. We have
  \[
  C_1^4 ( C_0, C_j ) = \sum_{l, k} n_{01}^1 n_{11}^l n_{l1}^{k} n_{k1}^{j}.
  \]
  From Lemma \ref{support 3}, for each $l \neq 0$ then exists $k \neq 0$ such that 
  \[n_{l1}^{k}n_{k1}^{j}> 0. \]   
  But  $n_{u1}^v > 0$ then $n_{u1}^{v} \geqslant
  1/(q + 1)$. Moreover
  \[\Pr ( l = 0 ) = n_{11}^0 = 1/(q + 1).\]
  Hence  
  \[  C_1^4 ( C_0, C_j ) \geqslant \frac{1}{( q + 1 )^2}
  \Pr ( l \neq 0 ) = \frac{q}{( q + 1 )^3}.\]  
   Therefore, we have
  \[C_1^4 ( C_0, C_j ) \geqslant \frac{q^3}{( q + 1 )^4}
  \frac{( q + 1 )}{q^2} > \frac{q^2 ( q - 1 )}{( q + 1 )^4} \pi ( C_j ).\]
  
 \item Suppose that $i \neq 0, j = 0$. Similar as in 2), we have 
  \begin{align*}C_1^4 ( C_i, C_0 ) &\geqslant \frac{q^3}{( q + 1 )^4}
  \frac{( q + 1 )}{q^2}\\ &> \frac{q^3}{( q + 1 )^4} \frac{1}{q^2} > \frac{q^2 (
  q - 1 )}{( q + 1 )^4} \pi ( C_j ).\end{align*}
  
 \item Suppose that $i, j \neq 0$. We have 
  \[ C_1^4 ( C_i, C_j ) = \sum_{t, l, k} n_{i1}^{t} n_{t1}^{l} n_{l1}^{k} n_{k1}^{j}. \]
  Similar as in 2), we have 
  \begin{align*}C_1^4 (C_i, C_j ) &\geqslant \frac{1}{( q + 1 )^2} \Pr ( l \neq 0 )\\ &= \frac{1}{( q +
  1 )^2} ( 1 - \Pr ( l = 0 ) ) .\end{align*} 
  But \[\Pr ( l = 0 ) = \sum_t n_{i1}^{t} n_{t1}^{0} = n_{i1}^{1} \leqslant \frac{2}{q + 1}.\]
  So we have \[C_1^4 ( C_i, C_j ) \geqslant \frac{q - 1}{( q
  + 1 )^3} = \frac{q^2 ( q - 1 )}{( q + 1 )^4} \pi ( C_j ).\]  
  \end{enumerate}
  This concludes the proof of the claim.
\end{proof}

From Lemma \ref{important 3}, we can determine the stationary distribution and the rate of convergence of the random walk $C_1$.

\begin{Theorem}\label{rate 3}
  Let $q \equiv 3$ (mod $4$). Then
  \[ \lim_{n \rightarrow \infty} C_1^n = \frac{1}{q^2} C_0 + \frac{q + 1}{q^2} \sum_{i \in F_q^\ast}C_i .\] 
  Furthermore, the rate of convergence (i.e. the mixing time of the random walk) is linear with respect to $q$. 
\end{Theorem}

\begin{proof}
  We create two copies of random walk $C_1$. The first one starts from $C_i$
  for fixed $i$ and the second one starts randomly in hypergroup of weighted-circles $C$
  with distribution $\pi$ in Lemma \ref{important 3}. In step $m^{\text{th}}$, suppose that we are in the 
  weighted-circle $C_t$ in the first walk and in the weighted-circle $C_s$ in the second walk for some $s$ and $t$.
  If $t = s$ then in the next step, we choose the step of the second walk which is the same with the first's. 
  Otherwise, let they walk by $C_1$ independently. It is clearly that both random walks have the same Markov kernel       $C_1$ and the second one has the distribution $\pi$.
    
  Set $c = 1 - q^2 ( q - 1 )/( q+1 )^4$. From Lemma \ref{important 3}, we have
  \begin{align*}  
  d_{\tmop{TV}} ( C_1^4 ( &C_i, . ), \pi ) = \frac{1}{2} \sum_j |
  \pi ( C_j ) - C_1^4 ( C_i, C_j ) | \\
  &= \sum_{j : C_1^4 ( C_i, C_j ) < \pi ( C_j
  )} ( \pi ( C_j ) - C_1^4 ( C_i, C_j ) )\\  
  &\leqslant \sum_{j : C_1^4 ( C_i, C_j ) < \pi
  ( C_j )} \pi ( C_j ) ( 1 - c ) \\
  &\leqslant 1 - c.
  \end{align*}
  Applying Lemma \ref{lap} we have
  \[d_{\tmop{TV}} ( C_1^{4 n} ( C_i, . ), \pi ) \leqslant ( 1 - c)^n.\]
  Thus, if $( 1 - c )^n < 1/2e$ then  
  \[d_{\tmop{TV}} ( C_1^{4 n} ( C_i, . ), \pi ) \leqslant 1/2e\]
  and $\tau_q \leqslant 4 n$.
  But, the inequality \[(1-c)^n = ( 1 - q^2 ( q - 1 )/( q+1 )^4)^n <
  1/2 e\] is equivalent to
  \[n \log \left( 1 - \frac{q^2 ( q - 1 )}{( q+1)^4} \right) < - \ln 2 - 1. \]
  This implies that
 \[n \left( \frac{q^2 ( q - 1 )}{( q + 1 )^4} + \frac{q^4 (
  q - 1 )^2}{( q+1 )^8} + \ldots \right) > 1 + \ln 2.\]
  Thus, we can choose 
 \[n > \frac{( 1 + \log 2 ) ( q + 1 )^4}{q^2 ( q - 1 )}.\]
 This conludes the proof of the theorem.
\end{proof}

Note that $|C_{0}|=1,|C_{i}|=q+1$ for $i\in F_q^\ast$, and the space $F_q^2$ has
$q^{2}$ points, so the distribution of $C_1^{n}$ is, in some sense, close to
uniform over the space $F_{q}^{2}$ when $n$ tends to infinite. Walking
randomly by any $C_{i}$ with $i\in F_q^\ast$ we have the same results as for $C_{1}$. In
hypergroup language the limiting distribution is the Haar measure on the hypergroup.

\subsubsection{The second case}
Suppose that $q \equiv 1$ (mod $4$). We have the following lemma.

\begin{Lemma} \label{important 1} Let $q \equiv 1$ (mod $4$) $\geqslant 13$ and $\pi$ be the distribution on the set $C = \{ C_i\}_{i \in F_q^{+}}$ with $\pi(C_0)=1/q^2$, $\pi(C_q)=(2q-1)/q^2$ and $\pi(C_i) = (q+1)/q^2$ for all $i \in F_q^\ast$. Then
  \[C_1^6 ( C_i, C_j ) \geqslant \frac{1}{3 q} \pi ( C_j )\]
  for all $i, j \in F_q^{+}$.
\end{Lemma}

\begin{proof}
  We have
  \begin{align*} C_1^6 ( C_i, C_j ) &= \sum_{k, l, m, t, h \in F_q^{+}} n_{i
     1}^k n_{k 1}^l n_{l 1}^m n_{m 1}^t n_{t 1}^h n_{h 1}^j\\ &\geqslant
     \sum_{k, l, m, t, h \in F_q^{+}, k, h \neq 0} n_{i 1}^k n_{k
     1}^l n_{l 1}^m n_{m 1}^t n_{t 1}^h n_{h 1}^j.\end{align*}
  For fixed $k, h \neq 0$, we want to approximate
  \[ \sum_{l, m, t \in F_q^{+}} n_{k 1}^l n_{l 1}^m n_{m 1}^t
     n_{t 1}^h. \]
  From Theorem \ref{hyper 1} , $n_{k 1}^l \leqslant \frac{2}{q - 1}$ for all $l$ so
  \[ \Pr ( l \neq 0, 1, q ) = \sum_{l \neq 0, 1, q} n_{k 1}^l \geqslant 1 -
     \frac{6}{q - 1} . \]
  Similarly, we also have Pr$( t \neq 0, 1, q ) \geqslant 1 - \frac{6}{q -
  1}$. We fix $l, t \neq 1, 0, q$. By Theorem \ref{hyper 1}, we have
  \[n_{l 1}^q = n_{q 1}^t = \frac{2}{q - 1}.\]
  Hence
  \begin{align*} \sum_{l, m, t \in F_p^{+}} n_{k 1}^l n_{l 1}^m n_{m 1}^t n_{t
     1}^h & \geqslant \sum_{l, t \neq 0, 1, q} n_{k 1}^l n_{l 1}^q n_{q 1}^t
     n_{t 1}^h\\ &\geqslant \frac{4}{( q - 1 )^2} \left( 1 - \frac{6}{q - 1}
     \right)^2. \end{align*}
  But $q \geqslant 13$ so
  \[ \frac{4}{( q - 1 )^2} \left( 1 - \frac{6}{q - 1} \right)^2 \geqslant \frac{1}{( q - 1 )^2}. \]  
  Therefore, if $k, h \neq 0$ then
  \[ \sum_{l, m, t \in F_q^{+}} n_{k 1}^l n_{l 1}^m n_{m 1}^t n_{t
     1}^h \geqslant \frac{1}{( q - 1 )^2}. \]
  This implies that
  \begin{align*} C_1^6 ( C_i, C_j ) &= \sum_{k, l, m, t, h \in F_q^{+}} n_{i
     1}^k n_{k 1}^l n_{l 1}^m n_{m 1}^t n_{t 1}^h n_{h 1}^j\\ &\geqslant
     \sum_{k, h \in F_q^{+}, \ne 0} n_{i 1}^k \times \frac{1}{( q - 1 )^2} \times
     n_{h 1}^j. \end{align*}
  Hence
  \[ C_1^6 ( C_i, C_j ) \geqslant \frac{1}{( q - 1 )^2} \Pr ( k \neq 0 ) \Pr
     ( h \neq 0 ). \]
  But from Theorem \ref{hyper 1}, $n_{i 1}^0, n_{01}^j \leqslant \frac{2}{q - 1}$. Thus, we have
  \[\Pr
  ( k \neq 0 ), \Pr ( h \neq 0 ) \geqslant 1 - \frac{2}{q - 1} \geqslant 1 -
  \frac{2}{12} = \frac{5}{6}.\]
  This implies that 
  \begin{align*}C_1^6 ( C_i, C_j ) &\geqslant
  \frac{25}{36 ( q - 1 )^2} = \frac{25 q^2}{72 ( q - 1 )^3} \frac{2 ( q
  - 1 )}{q^2}\\ &\geqslant \frac{25 q^2}{72 ( q - 1 )^3} \pi ( C_j ) > \frac{1}{3
  q} \pi ( C_j ).\end{align*}
  This concludes the proof of the lemma.
\end{proof}

Similarly, from Lemma \ref{important 1}, we can determine the stationary distribution and the rate of convergence of the random walk $C_1$. 

\begin{Theorem}\label{rate 1}
  Let $q \equiv 1$ (mod $4$). Then
  \[ \lim_{n \rightarrow \infty} C_1^n = \frac{1}{q^2} C_0 + \frac{2(q-1)}{q^2} C_q + \frac{q - 1}{q^2} \sum_{i \in F_q^\ast}C_i .\] 
  Furthermore, the rate of convergence (i.e. the mixing time of the random walk) is linear with respect to $q$. 
\end{Theorem}

The proof of this theorem is omitted since it is the same as the proof of Theorem \ref{rate 3}. Note that,
walking randomly by any $C_{i}$ with $i\in F_q^\ast$ we have the same results as for $C_{1}$ and the limiting distribution in Theorem \ref{rate 1} is the Haar measure on the hypergroup.

\end{document}